\documentclass[11pt]{amsart}
\usepackage{bbm}
\usepackage{amssymb, amsthm, amsmath, amscd}
\usepackage{enumerate}
\usepackage[usenames,dvipsnames]{color}
\usepackage{hyperref}
\hypersetup{hypertex=true,
	colorlinks=true,
	linkcolor=blue,
	anchorcolor=blue,
	citecolor=blue}

\numberwithin{equation}{section}

\theoremstyle{cupthm}
\newtheorem{thm}{Theorem}[section]
\newtheorem{prop}[thm]{Proposition}
\newtheorem{cor}[thm]{Corollary}
\newtheorem{lemma}[thm]{Lemma}
\theoremstyle{cupdefn}
\newtheorem{defn}[thm]{Definition}
\theoremstyle{cuprem}

\numberwithin{equation}{section}

\usepackage{amsfonts}
\usepackage{mathrsfs}
\usepackage{textcomp}
\usepackage[all]{xy}

\title[ Uniform property $\Gamma$ for Crossed products]{\bf Uniform property $\Gamma$ for Crossed products by group actions with the Rokhlin-type properties }

\author{Xiaochun Fang}
\address{School of Mathematical Sciences,
	Key Laboratory of Intelligent Computing and Applications(Ministry of Education), Tongji University, Shanghai 200092}
\email{xfang@tongji.edu.cn}

\author{Haotian Tian}
\address{School of Mathematical Sciences,
	Key Laboratory of Intelligent Computing and Applications(Ministry of Education), Tongji University, Shanghai 200092}
\email{2011284@tongji.edu.cn}


\begin{document}
	\begin{abstract}
		In this paper, let $A$ be a unital separable simple infinite dimensional C*-algebra which has uniform property $\Gamma$. Let $\alpha\colon G\to \mathrm{Aut}(A)$ be an action of a finite group which has the weak tracial Rokhlin property. Then we prove that the crossed product $A\rtimes_\alpha G$ and fixed point algebra $A^\alpha$ have uniform property $\Gamma$. Let $\alpha\colon G\to \mathrm{Aut}(A)$ be an action of a second-countable compact group which has the tracial Rokhlin property with comparison. Then we prove that the crossed product $A\rtimes_\alpha G$ and fixed point algebra $A^\alpha$ have uniform property $\Gamma$.
	\end{abstract}
	\maketitle
	\section{Introduction}
	The Rokhlin property for the case of a single automorphism was originally introduced for von Neumann algebras by Connes in \cite{AC1975}. Later, the Rokhlin property for finite group actions on C*-algebras first appeared in the work of Herman and Jones in \cite{HJ1982} and \cite{HJ1983}. This property is useful to understand the structure of the crossed product of C*-algebras and properties passing from the original algebra to the crossed product \cite{OP2012}. However, the finite group acitions with the Rokhlin property are rare. Phillips, in \cite{NCP2011}, introduced the tracial Rokhlin property for finite group actions on unital simple C*-algebras. The tracial Rokhlin property is generic in many cases, and also can be used to study properties passing from the original algebra to the crossed product. Weak versions of the tracial Rokhlin property in which one uses orthogonal positive contractions instead of orthogonal projections were studied for actions on unital simple C*-algebras with few projections \cite{NCP2012,MS2012,HO2013,QW2013,GHS2021,QW2018} (see Definition \ref{WTRF}). As an example, the flip action on the Jiang-Su algebra $\mathcal{Z}\cong \mathcal{Z}\otimes\mathcal{Z}$ has the weak tracial Rokhlin property but it does not have the tracial Rokhlin property \cite{HO2013}. For the non-unital case, Santiago and Gardella studyed the Rokhlin property for finite group actions on non-unital simple C*-algebras in \cite{LS2015} and \cite{GS2016}. Forough and Golestani studyed the (weak) tracial Rokhlin property for finite group actions on non-unital simple C*-algebras in \cite{FG2020}.
	
	In \cite{HW2007}, Hirshberg and Winter also introduced the Rokhlin property for second-countable compact group actions on unital C*-algebras. Since then, crossed products by compact group actions with the Rokhlin property have been studied by several authors. In particular, permanence properties are proved in \cite{HW2007}, \cite{EG2017} and \cite{EG2019}. The same as finite groups, Rokhlin actions of compact groups are rare, especially when the group is connected. More recently, Mohammadkarimi and Phillips studied the tracial Rokhlin property with comparison for compact group actions and proved that the crossed product of a unital separable simple infinite dimensional C*-algebra with tracial rank zero by an action of a second-countable compact group with the tracial Rokhlin property with comparison has again tracial rank zero in \cite{MP2022} and some other permanence properties. Moreover, they gave some examples of compact group actions with the tracial Rokhlin property with comparison. The authors have studied many permanence properties in \cite{TF2022} including stable rank one, real rank zero, $\beta$-comparison, Winter's $n$-comparison, $m$-almost divisibility and weakly ($m$,$n$)-divisibility.
	
	The Elliott program aims to classify amenable C*-algebras. In his efforts to classify simple separable amenable C*-algebras, Elliott highlighted the necessity of considering certain regularity properties of these algebras. Three particular properties of interest are: finite nuclear dimension, tensorial absorption of the Jiang–Su algebra $\mathcal{Z}$ (also known as $\mathcal{Z}$-stability), and strict comparison of positive elements. Toms and Winter conjectured, in what is known as the Toms–Winter conjecture (see, e.g., \cite{ET2008}), that these three fundamental properties are equivalent for all separable, simple, unital, amenable C*-algebras. This conjecture has now been nearly entirely proven (see \cite{BBSTWW2019,CETWW2021,KR2014,TWW2015}).
	
	To prove that $\mathcal{Z}$-stability implies finite nuclear dimension, Castillejos et al. introduced the uniform property $\Gamma$ and the complemented tracial orthogonal partitions of unity property for separable C*-algebras in \cite{CETWW2021}. They showed that $\mathcal{Z}$-stability implies the uniform property $\Gamma$, and that the uniform property $\Gamma$ in turn implies the complemented tracial orthogonal partitions of unity property and this can prove finite nuclear dimension for separable simple nuclear nonelementary unital C*-algebra. Also, in \cite{CETWW2021}, they showed that the Toms–Winter conjecture holds for separable simple unital non-elementary C*-algebras that have the uniform property $\Gamma$.
	
	Examples of separable amenable C*-algebras with the uniform property $\Gamma$ are now abundant. Kerr and Szab\'o established the uniform property $\Gamma$ for crossed product C*-algebras that arise from a free action of an infinite amenable group with the small boundary property on a compact metrizable space (see \cite[Theorem 9.4]{KS2020}).
	
	In this paper, we get the following results:
	
	\begin{thm}
		Let $A$ be a unital separable simple infinite dimensional C*-algebra which has uniform property $\Gamma$. Let $\alpha\colon G\to \mathrm{Aut}(A)$ be an action of a finite group which has the weak tracial Rokhlin property. Then the crossed product $A\rtimes_\alpha G$ and the fixed point algebra $A^\alpha$ have uniform property $\Gamma$.
	\end{thm}
	
	\begin{thm}
		Let $A$ be a unital separable simple infinite dimensional C*-algebra which has uniform property $\Gamma$. Let $\alpha\colon G\to \mathrm{Aut}(A)$ be an action of a second-countable compact group which has the tracial Rokhlin property with comparison. Then the crossed product $A\rtimes_\alpha G$ and the fixed point algebra $A^\alpha$ have uniform property $\Gamma$.
	\end{thm}
	
	The paper is organized as follows. Section 2 contains some preliminaries about ultraproducts, limit traces, Cuntz subequivalence and actions with the Rokhlin-type properties. Section 3 contains the proofs of the main theorems and corollarys.
	
	
	\section{Preliminaries and Definitions}
	In this section, we recall some definitions and known facts about ultraproducts, limit traces, Cuntz subequivalence and actions with the Rokhlin-type properties.
	
	\begin{defn}
		Let $\omega\in\beta\mathbb{N}\backslash \mathbb{N}$ be a fixed free ultrafilter. Let $A$ be a C*-algebra. We use $l^\infty(\mathbb{N},A)$ to denote the set of all bounded sequences in $A$ with the supremum norm. The ultrapower of $A$ is then given by
		\[A_\omega=l^\infty(\mathbb{N},A)/\{(a_n)_{n\in\mathbb{N}}\colon\lim_{n\to\omega}\|a_n\|=0\}.\]
		Denoted by $\kappa_A\colon l^\infty(\mathbb{N},A)\to A_\omega$ the quotient map. Define $\iota\colon A\to l^\infty(\mathbb{N},A)$ by $\iota(a)=(a,a,a,\dots)$, the constant sequence, for all $a\in A$. Identify $A$ with $\kappa_A\circ\iota(A)$. We will adopt a standard abuse of notation and denote elements in $A_\omega$ by choice of a representative sequence $(a_n)$. 
		
		A tracial state $\tau$ on $A_\omega$ is called a limit trace if there is a sequence $(\tau_n)$ of tracial states on $A$ such that $\tau((a_n)_{n\in\mathbb{N}})=\lim_{n\to\omega}\tau_n(a_n)$ for all $(a_n)_{n\in\mathbb{N}}\in A_\omega$. The set of limit traces on $A_\omega$ will be denoted by $\mathrm{T}_\omega(A)$.
		
		Suppose that $\mathrm{T}(A)$ is non-empty, the trace kernel ideal is given by
		\[J_A=\{x\in A_\omega\colon\tau(x^*x)=0 \ \mathrm{for}\ \mathrm{all}\ \tau\in \mathrm{T}_\omega(A)\}.\]
		The uniform tracial ultrapower of $A$ is defined as
		\[A^\omega=A_\omega/J_A.\]
		When $A$ is separable, $A^\omega$ is unital if and only if $\mathrm{T}(A)$ is compact by \cite[Proposition
		1.11]{CETWW2021}. The notation $\mathrm{T}_\omega(A)$ will also be used for tracial states on $A^\omega$ coming from limit traces. There is a canonical map $\iota'\colon A\to A^\omega$ given by taking constant sequences. This need not be an embedding in general, but it will be whenever $\mathrm{T}(A)$ is separating. Abusing notation slightly, we will simply write $A^\omega\cap A'$ instead of $A^\omega\cap \iota'(A)'$.
	\end{defn}
	
	\begin{defn}
		Let $A$ be a C*-algebra, $a\in A_+$ and $\varepsilon>0$. Then we denote $f(a)$ by $(a-\varepsilon)_+$, where $f(t)=max\{0,t-\varepsilon\}$ is continuous from $[0,\infty)$ to $[0,\infty)$.
	\end{defn}
	
	The following definitions related to Cuntz comparison are from \cite{JC1978}, for more information, you can refer to \cite{GP2022} and \cite{APT2018}.
	
	\begin{defn}
		Let $A$ be a C*-algebra. Let $a,b\in (A\otimes K)_+$.
		\begin{enumerate}
			\item We say that $a$ is Cuntz subequivalent to $b$ (written $a\precsim_A b$), if there is a sequence $(r_n)_{n=1}^\infty$ in $A\otimes K$ such that $\lim\limits_{n\to\infty}\|r_n^*br_n-a\|=0$.
			
			\item We say that $a$ is Cuntz equivalent to $b$ (written $a\sim_A b$), if $a\precsim_A b$ and $b\precsim_A a$. This is an equivalence relation, we use $\langle a\rangle_A$ to denote the equivalence class of $a$. With the addition operation $\langle a\rangle_A +\langle b\rangle_A=\langle a\oplus b\rangle_A$ and the order operation $\langle a\rangle_A \leq\langle b\rangle_A$ if $a\precsim_A b$, $\mathrm{Cu}(A)=(A\otimes K)_+/\sim_A$ is an ordered semigroup which we called Cuntz semigroup. $\mathrm{W}(A)=M_\infty(A)_+/\sim_A$ is also an ordered semigroup with the same operation and order as above.
		\end{enumerate}
		If $B$ is a hereditary C*-subalgebra of $A$, and $a,b\in B_+$, then it is easy to check that $a\precsim_A b\iff a\precsim_B b$.
	\end{defn}
	
	\begin{defn}\cite[Definition 1.3]{HP2015}
		Let $G$ be a compact group, and let $A$ be a C*-algebra, $B$ be a C*-algebra. Let $\alpha\colon G\to \mathrm{Aut}(A)$ and $\gamma\colon G\to \mathrm{Aut}(B)$ be actions of $G$ on $A$ and $B$. Let $F\subseteq A$ and $S\subseteq B$ be subsets, and let $\varepsilon>0$. A completely positive contractive map $\varphi\colon A\to B$ is said to be an ($F$,$S$,$\varepsilon$)-approximately central equivariant multiplicative map if:
		\begin{enumerate}
			\item $\|\varphi(xy)-\varphi(x)\varphi(y)\|<\varepsilon$ for all $x,y\in F$.
			
			\item $\|\varphi(x)a-a\varphi(x)\|<\varepsilon$ for all $x\in F$ and all $a\in S$.
			
			\item $\sup_{g\in G}\|\varphi(\alpha_g(x))-\gamma_g(\varphi(x))\|<\varepsilon$ for all $x\in F$.
		\end{enumerate}
	\end{defn}
	
	\begin{defn}\cite[Definition 1.4]{MP2022}
		Let $A$ and $B$ be C*-algebras, and let $F\subseteq A$. A completely positive contractive map $\varphi\colon A\to B$ is said to be an ($n$,$F$,$\varepsilon$)-approximately multiplicative map if whenever $m\in\{1,2,\dots,n\}$ and $x_1,x_2,\dots,x_m\in F$, we have
		\[\|\varphi(x_1x_2\dots x_m)-\varphi(x_1)\varphi(x_2)\dots\varphi(x_m)\|<\varepsilon.\]
		If $S\subseteq B$ is also given, then $\varphi$ is said to be an ($n$,$F$,$S$,$\varepsilon$)-approximately central multiplicative map if, in addition, $\|\varphi(x)a-a\varphi(x)\|<\varepsilon$ for all $x\in F$ and all $a\in S$.
	\end{defn}
	
	Now, let us recall the notion of the tracial Rokhlin property with comparison for second-countable compact group actions defined by Mohammadkarimi and Phillips in \cite{MP2022}.
	
	\begin{defn}\cite[Definition 2.4]{MP2022}\label{DEF}
		Let $G$ be a second-countable compact group, let $A$ be a unital simple infinite dimensional C*-algebra, and let $\alpha\colon G\to \mathrm{Aut}(A)$ be an action. We say that the action $\alpha$ has the tracial Rokhlin property with comparison if for any $\varepsilon>0$,any finite set $F\subseteq A$, any finite set $S\subseteq C(G)$, any $x\in A_+$ with $\|x\|=1$, and any $y\in(A^\alpha)_+\setminus\{0\}$, there exist a projection $p\in A^\alpha$ and a unital completely positive map $\psi\colon C(G)\to pAp$ such that
		\begin{enumerate}
			\item $\psi$ is an ($F$,$S$,$\varepsilon$)-approximately central equivariant multiplicative map.
			
			\item $1-p\precsim_A x$.
			
			\item $1-p\precsim_{A^\alpha} y$.
			
			\item $1-p\precsim_{A^\alpha} p$.
			
			\item $\|pxp\|>1-\varepsilon$.
		\end{enumerate}
	\end{defn}
	
	The next theorem is the key tool for transferring properties from the original algebra to the fixed point algebra.
	\begin{thm}\cite[Theorem 2.17]{MP2022}\label{MTF}
		Let $G$ be a second-countable compact group, let $A$ be a unital separable simple infinite dimensional C*-algebra and let $\alpha\colon G\to \mathrm{Aut}(A)$ be an action with the tracial Rokhlin property. Then for any $\varepsilon>0$, any $n\in \mathbb{N}$, any compact subset $F_1\subseteq A$, any compact subset $F_2\subseteq A^\alpha$, any $x\in A_+$ with $\|x\|=1$, and any $y\in (A^\alpha)_+\setminus\{0\}$, there exist a projection $p\in A^\alpha$ and a unital completely positive map $\varphi\colon A\to pA^\alpha p$ such that
		\begin{enumerate}
			\item $\varphi$ is an ($n$,$F_1\cup F_2$,$\varepsilon$)-approximately multiplicative map.
			
			\item $\|pa-ap\|<\varepsilon$ for all $a\in F_1\cup F_2$.
			
			\item $\|\varphi(a)-pap\|<\varepsilon$ for all $a\in F_2$.
			
			\item $\|\varphi(a)\|\geq\|a\|-\varepsilon$ for all $a\in F_1\cup F_2$.
			
			\item $1-p\precsim_A x$.
			
			\item $1-p\precsim_{A^\alpha} y$.
			
			\item $1-p\precsim_{A^\alpha} p$
			
			\item $\|pxp\|>1-\varepsilon$.
		\end{enumerate}
	\end{thm}
	
	\begin{thm}\cite[Theorem 3.9, Corollary 3.10]{MP2022}\label{CPS}
		Let $G$ be a second-countable compact group, let $A$ be a unital separable simple infinite dimensional C*-algebra and let $\alpha\colon G\to \mathrm{Aut}(A)$ be an action with the tracial Rokhlin property with comparison. Then the crossed product $A\rtimes_\alpha G$ is simple. Moreover, the algebras $A\rtimes_\alpha G$ and $A^\alpha$ are Morita equivalent and stably isomorphic.
	\end{thm}
	
	The notion of the weak tracial Rokhlin property with comparison for finite group actions was introduced by Asadi-Vasfi, Golestani, Phillips in \cite{ANP2021}. 
	\begin{defn}\cite[Definition 3.2]{ANP2021}\label{WTRF}
		Let $G$ be a finite group, let $A$ be a unital simple infinite dimensional C*-algebra, and let $\alpha\colon G\to \mathrm{Aut}(A)$ be an action. We say that $\alpha$ has the weak tracial Rokhlin property if for any finite set $F\subseteq A$, any $\varepsilon>0$, any $x\in A_+$ with $\|x\|=1$, there exist orthogonal positive contractions $(d_g)_{g\in G}\in A$ with $d=\sum_{g\in G}d_g$ such that
		\begin{enumerate}
			\item$\|d_ga-ad_g\|<\varepsilon$ for all $a\in F$ and all $g\in G$.
			
			\item$\|\alpha_g(d_h)-d_{gh}\|<\varepsilon$ for all $g,h\in G$.
			
			\item$1-d\precsim_A x$.
			
			\item$\|dxd\|>1-\varepsilon$.
		\end{enumerate}
	\end{defn}
	
	The next theorem is an approximation property (and in a slightly different form in the earlier \cite[Lemma VII.4, Lemma VII.16]{DA2008}) which is closely related to the notion of essential tracial approximation (see \cite[Definition 3.1]{FL2022}) and the notion of generalized tracial approximation (see \cite[Definition 1.2]{EFF2023}).
	\begin{thm}\cite[Theorem 3.4]{FW2024}\label{GTA1}
		Let $A$ be a unital simple infinite dimensional C*-algebra. Let $\alpha\colon G\to \mathrm{Aut}(A)$ be an action of finite group which has the weak tracial Rokhlin property. Then for any finite subset $F\subset A\rtimes_\alpha G$, any $\varepsilon>0$, and any nonzero positive $x\in A\rtimes_\alpha G$, there exist a positive contraction $f\in A$, a C*-subalgebra $B$ of $A\rtimes_\alpha G$ with $B\cong \overline{fAf}\otimes M_l$ ($l$=Card($G$)) and a positive contraction $d\in B$ such that
		\begin{enumerate}
			\item $\|da-ad\|<\varepsilon$ for all $a\in F$.
			\item $da\in_\varepsilon B$ for all $a\in F$.
			\item $1-d\precsim_{A\rtimes_\alpha G} x$.
			\item $\|dad\|>\|a\|-\varepsilon$ for all $a\in F$. 
		\end{enumerate}
	\end{thm}
	
	Like in \cite[Lemma 3.5]{FW2025}, we can do some functional calculus for $d$ in Theorem \ref{GTA1}.
	
	\begin{lemma}($\mathrm{cf}$.\cite[Lemma 3.5]{FW2025}\label{GTA}
		Let $A$ be a unital simple infinite dimensional C*-algebra. Let $\alpha\colon G\to \mathrm{Aut}(A)$ be an action of finite group which has the weak tracial Rokhlin property. Let $f,g\colon [0,1]\to \mathbb{C}$ be continuous functions with $g(0)=0$. Then for any finite subset $F\subset A\rtimes_\alpha G$, any $\varepsilon>0$, and any nonzero positive $x\in A\rtimes_\alpha G$, there exist a positive contraction $f\in A$, a C*-subalgebra $B$ of $A\rtimes_\alpha G$ with $B\cong \overline{fAf}\otimes M_l$ ($l$=Card($G$)) and a positive contraction $d\in B$ such that
		\begin{enumerate}
			\item $\|f(d)a-af(d)\|<\varepsilon$ for all $a\in F$.
			\item $g(d)a\in_\varepsilon B$ for all $a\in F$. Moreover, if $a\in F$ is positive, then there exists a positive element $b\in B$ such that $\|g(d)ag(d)-b\|<\varepsilon$.
			\item $1-d\precsim_{A\rtimes_\alpha G} x$.
			\item $\|dad\|>\|a\|-\varepsilon$ for all $a\in F$. 
		\end{enumerate}
		\begin{proof}
			The proof is the same as that of \cite[Lemma 3.5]{FW2025}, so we omit it.
		\end{proof}
	\end{lemma}
	
	Uniform property $\Gamma$ was introduced by Castillejos et al., that was used to prove that $\mathcal{Z}$-stable implies that finite nuclear dimension in \cite{CETWW2021}.
	\begin{defn}\cite[Definition 2.1]{CETWW2021}
		Let $A$ be a separable C*-algebra with $\mathrm{T}(A)$ non-empty and compact. Then $A$ is said to have uniform property $\Gamma$ if for all $n\in\mathbb{N}$, there exist projections $p_1,\dots,p_n\in A^\omega\cap A'$ summing to $1_{A_\omega}$, such that
		\[\tau(ap_i)=\frac{1}{n}\tau(a),\ a\in A,\ \tau\in\mathrm{T}_\omega(A),\ i=1,\dots,n.\]
	\end{defn}
	We recall that the equivalent local refinement of uniform property $\Gamma$ from \cite[Proposition 2.4]{CETW2022}.
	\begin{prop}\cite[Proposition 2.4]{CETW2022}\label{UPG}
		Let $A$ be a separable C*-algebra with $\mathrm{T}(A)$ non-empty and compact. Then the following are equivalent:
		\begin{enumerate}
			\item $A$ has uniform property $\Gamma$.
			\item For any finite subset $F\subset A$, any $\varepsilon>0$, and any $n\in\mathbb{N}$, there exist pairwise orthogonal positive contractions $e_1,\dots,e_n\in A$ such that for $i=1,\dots,n$ and $a\in F$, we have
			\[\|e_ia-ae_i\|<\varepsilon\ \mathrm{and}\ \sup_{\tau\in\mathrm{T}(A)}|\tau(ae_i)-\frac{1}{n}\tau(a)|<\varepsilon.\]
		\end{enumerate}
	\end{prop}
	
	The next lemma will be used several times, so we state it here.
	\begin{lemma}\cite[Lemma 2.5.12]{HL2001}\label{RO}
		For any $\varepsilon>0$ and any integer $n>0$, there exists $\delta(\varepsilon,n)>0$ satisfying the following: If $A$ is a C*-algebra and $a_1,\cdots,a_n\in A_+$ with $\|a_i\|\leq 1 (i=1,\cdots,n)$ such that $\|a_ia_j\|<\delta$ when $i\neq j$, then there are $b_1,\cdots,b_n\in A_+$ such that $b_ib_j=0$ when $i\neq j$ and $\|a_i-b_i\|<\varepsilon$, $i=1,\cdots,n$.
	\end{lemma}
		
	\section{The main results}
	In this section, we give the proof of the main theorem. Before that, we give some basic propositions and lemmas.
	
	
	The following lemma is trival but we still state it here.
	\begin{lemma}($\mathrm{cf}$.\cite[Proposition 3.4]{HL2023})\label{TM}
		Let $A$ be a separable C*-algebra with $\mathrm{T}(A)$ nonempty and compact. Suppose that $A$ has uniform property $\Gamma$. Then, for any $k\in\mathbb{N}$, $M_k(A)$ also has uniform property $\Gamma$.
		\begin{proof}
			By using $\mathrm{T}(A)$ instead of $\mathrm{QT}(A)$, the proof is the same as that of \cite[Proposition 3.4]{HL2023}.
		\end{proof}
	\end{lemma}
	Although uniform property $\Gamma$ does not pass to hereditary C*-subalgebras in general, we will show that it holds when $A$ is unital simple and separable.
	\begin{lemma}\label{HS}
		Let $A$ be a unital simple separable C*-algebra with $\mathrm{T}(A)$ nonempty and compact. Suppose that $A$ has uniform property $\Gamma$. Then, for any hereditary C*-subalgebra $B$ of $A$, $B$ also has uniform property $\Gamma$.
		\begin{proof}
			Since $A$ is separable, we know that both $A$ and $B$ are $\sigma$-unital. Since $A$ is simple, we know that $B$ is a full hereditary C*-subalgebra of $A$. Thus, they are stably isomorphic. It follows from \cite[Proposition 2.6]{CE2021} that $A$ has stablized property $\Gamma$. By \cite[Theorem 2.10]{CE2021}, we know that $B$ has stablized property $\Gamma$. It follows from \cite[Proposition 2.6]{CE2021} that $B$ has uniform property $\Gamma$.
		\end{proof}
	\end{lemma}
	
	Now, we will give the proof of the first theorem of our main results. We mainly use the method of tracial approximation and consider the relation of tracial states spaces between the crossed product and its C*-subalgebra.
	\begin{thm}\label{FG}
		Let $A$ be a unital separable simple infinite dimensional C*-algebra which has uniform property $\Gamma$. Let $\alpha\colon G\to \mathrm{Aut}(A)$ be an action of a finite group which has the weak tracial Rokhlin property. Then then crossed product $A\rtimes_\alpha G$ has uniform property $\Gamma$.
		\begin{proof}
			Since $A$ has uniform property $\Gamma$, we know that  $\mathrm{T}(A)$ is nonempty. For $\tau\in\mathrm{T}(A)$, we can restrict it on $A^\alpha$ as a trace. Since $1_A=1_{A^\alpha}$, we know that $\tau\in\mathrm{T}(A^\alpha)$. Thus $\mathrm{T}(A^\alpha)$ is nonempty. By \cite[Proposition 3.2]{FG2020}, $\alpha$ is point wise outer, and so $A\rtimes_\alpha G$ is simple. Thus, $A^\alpha$ is Morita equivalent to $A\rtimes_\alpha G$. Since both algebras are separable and unital, $A\rtimes_\alpha G\cong pM_m(A^\alpha)p$ for some $m\in\mathbb{N}$ and $p\in M_m(A^\alpha)$. Therefore, $\mathrm{T}(A\rtimes_\alpha G)$ is nonempty. This together with the unitality of $A\rtimes_\alpha G$ implies that $\mathrm{T}(A\rtimes_\alpha G)$ is compact. By Proposition \ref{UPG}, we need to show that for any finite subset $F=\{a_1,\dots,a_k\}$ of $A\rtimes_\alpha G$, any $\varepsilon>0$, and any $n\in\mathbb{N}$, there exist pairwise orthogonal positive contractions $e_1,\dots,e_n\in A\rtimes_\alpha G$ such that
			\[\|e_ia_j-a_je_i\|<\varepsilon \quad \mathrm{and} \quad \sup_{\tau\in\mathrm{T}(A\rtimes_\alpha G)}|\tau(a_je_i)-\frac{1}{n}\tau(a_j)|<\varepsilon,
			\]
			for $i=1,\dots,n$ and $j=1,\dots,k$.
			
			Without loss of generality, we may assume that $\|a\|\leq1$ for all $a\in F$. We choose $\delta=\min(\frac{\varepsilon}{18},\frac{n\varepsilon}{12})$. Since $A\rtimes_\alpha G$ is simple and not of type \uppercase\expandafter{\romannumeral1}, by \cite[Corollary 2.5]{NCP2016}, there exists a nonzero positive element $x\in A\rtimes_\alpha G$ such that $d_\tau(x)<\delta$ for all $\tau\in\mathrm{T}(A\rtimes_\alpha G)$.
			
			Apply Lemma \ref{GTA} for $F$, $\delta$ and $x$, we get an element $f\in A_+^1$, a C*-subalgebra $B$ of $A\rtimes_\alpha G$ with $B\cong \overline{fAf}\otimes M_l$ ($l$=Card($G$)) and an element $d\in B_+$ such that
			\begin{enumerate}
				\item $\|(1-d)^\frac{1}{2}a-a(1-d)^\frac{1}{2}\|<\delta$, $\|d^\frac{1}{2}a-ad^\frac{1}{2}\|<\delta$ and $\|da-ad\|<\delta$ for all $a\in F$.\label{item1.1}
				\item $da\in_\delta B$ for all $a\in F$.\label{item1.2}
				\item $1-d\precsim_{A\rtimes_\alpha G} x$.\label{item1.3}
				\item $\|dad\|>\|a\|-\delta$ for all $a\in F$. 
			\end{enumerate}
			
			By (\ref{item1.3}), we have 
			\begin{equation}\label{eq9}
				\tau(1-d)\leq d_\tau(1-d)\leq d_\tau(x)<\delta,
			\end{equation}
			for all $\tau\in \mathrm{T}(A\rtimes_\alpha G)$.
			
			By (\ref{item1.2}), there exist a finite subset $F'=\{a'_1,\dots,a'_k\}$ of $B$ such that
			\begin{equation}\label{eq7}
				\|da_j-a_j'\|<\delta,
			\end{equation}
			for all $j=1,\dots,k$. Put $a''_j=(1-d)^\frac{1}{2}a_j(1-d)^\frac{1}{2}$, for all $j=1,\dots,k$. Then, by (\ref{item1.1}) and (\ref{item1.2}), we have
			\begin{align*}
				\|a_j-a_j'-a_j''\|\leq&{}\|da_j-a_j'\|+\|(1-d)a_j-a_j''\|\\
				<&{}\delta+\|(1-d)a_j-(1-d)^\frac{1}{2}a_j(1-d)^\frac{1}{2}\|\\
				<&{}\delta+\delta=2\delta<\frac{\varepsilon}{3},
			\end{align*}
			for $j=1,\dots,k$.
			
			For $\frac{\varepsilon}{9}>0$ and $n\in\mathbb{N}$, we choose $\delta'<\min(\frac{\varepsilon}{9},\frac{2\varepsilon}{3})$ sufficiently small such that satisfying Lemma \ref{RO}.
			
			Since $A$ has uniform property $\Gamma$ and $B\cong \overline{fAf}\otimes M_l$, by Lemma \ref{TM} and Lemma \ref{HS}, we know that $B$ has uniform property $\Gamma$. Apply Proposition \ref{UPG} for $n$, $F''=F'\cup\{d^\frac{1}{2}\}$ and $\frac{\delta'}{2}$ as given, we have pairwise orthogonal positive contractions $e_1',\dots,e_n'\in B$ such that
			\begin{equation}\label{eq11}
				\|e_i'a-ae_i'\|<\frac{\delta'}{2} \quad \mathrm{and} \quad \sup_{\tau\in\mathrm{T}(B)}|\tau(ae_i')-\frac{1}{n}\tau(a)|<\frac{\delta'}{2},
			\end{equation}
			for $i=1,\dots,n$ and all $a\in F''$.
			
			Since $e_i'e_j'=0$ for $i\neq j$, by (\ref{eq11}), we have 
			\begin{align*}
				\|d^\frac{1}{2}e_i'd^\frac{1}{2}d^\frac{1}{2}e_j'd^\frac{1}{2}\|\leq&{}\|d^\frac{1}{2}e_i'd^\frac{1}{2}d^\frac{1}{2}e_j'd^\frac{1}{2}-de_i'd^\frac{1}{2}e_j'd^\frac{1}{2}\|+\|de_i'd^\frac{1}{2}e_j'd^\frac{1}{2}-de_i'e_j'd\|\\<&{}\frac{\delta'}{2}+\frac{\delta'}{2}=\delta'.
			\end{align*}
			By Lemma \ref{RO}, there exist pairwise orthogonal positive contractions $e_1,\dots,e_n\in B_+$ such that $\|e_i-d^\frac{1}{2}e_i'd^\frac{1}{2}\|<\frac{\varepsilon}{9}$ for $i=1,\dots,n$.
			
			Thus, using (\ref{eq11}), (\ref{item1.1}), (\ref{eq7}), we have
			\begin{align*}
				\|a_je_i-a'_je_i'\|\leq&{}\|a_je_i-a_jd^\frac{1}{2}e_i'd^\frac{1}{2}\|+\|a_jd^\frac{1}{2}e_i'd^\frac{1}{2}-da_je_i'\|+\|da_je_i'-a_j'e_i'\|\\
				<&{}\frac{\varepsilon}{9}+\|a_jd^\frac{1}{2}e_i'd^\frac{1}{2}-a_jde_i'\|+\|a_jde_i'-da_je_i'\|+\|da_je_i'-a_j'e_i'\|\\
				<&{}\frac{\varepsilon}{9}+\frac{\delta'}{2}+\delta+\delta\\
				<&{}\frac{\varepsilon}{9}+\frac{\varepsilon}{18}+\frac{2\varepsilon}{18}<\frac{\varepsilon}{3}.
			\end{align*}
			With the same arguement, we have
			\[\|e_ia_j-e_i'a_j'\|<\frac{\varepsilon}{3}.\]
			Since $\|e_i'a_j'-a_j'e_i'\|<\frac{\varepsilon}{3}$, using (\ref{eq7}) we have
			\begin{align*}
				\|a_je_i-e_ia_j\|\leq\|a_je_i-a'_je_i'\|+\|e_i'a_j'-a_j'e_i'\|+\|e_ia_j-e_i'a_j'\|
				<\frac{\varepsilon}{3}+\frac{\varepsilon}{3}+\frac{\varepsilon}{3}=\varepsilon.
			\end{align*}
			Since $\|a_je_i-a'_je_i'\|<\frac{\varepsilon}{3}$, for all $\tau\in \mathrm{T}(A\rtimes_\alpha G)$, we have
			\begin{equation}\label{eq8}
				|\tau(a_je_i)-\tau(a_j'e_i')|<\frac{\varepsilon}{3}.
			\end{equation}
			Since $\|a_i-a_i'-a_i''\|<\frac{\varepsilon}{3}$, we have
			\[|\tau(a_j)-\tau(a_j')-\tau(a_j'')|<\frac{\varepsilon}{3}.\]
			Note that $a_j''=(1-d)^\frac{1}{2}a_j(1-d)^\frac{1}{2}$ and using (\ref{eq8}) at the second step, (\ref{eq9}) at the fourth step, we have
			\begin{align*}
				|\tau(a_je_i)-\frac{1}{n}\tau(a_j)|\leq&{}|\tau(a_je_i)-\tau(a_j'e_i')|+|\tau(a_j'e_i')-\frac{1}{n}\tau(a_j')|+|\frac{1}{n}\tau(a_j')-\frac{1}{n}\tau(a_j)|\\
				<&{}|\tau(a_j'e_i')-\frac{1}{n}\tau(a_j')|+\frac{\varepsilon}{3}+\frac{1}{n}(|\tau(a_j'')|+3\delta)\\
				\leq&{}|\tau(a_j'e_i')-\frac{1}{n}\tau(a_j')|+\frac{\varepsilon}{3}+\frac{1}{n}(|\tau(1-d)|+3\delta)\\
				<&{}|\tau(a_j'e_i')-\frac{1}{n}\tau(a_j')|+\frac{\varepsilon}{3}+\frac{1}{n}4\delta\\
				\leq&{}|\tau(a_j'e_i')-\frac{1}{n}\tau(a_j')|+\frac{\varepsilon}{3}+\frac{\varepsilon}{3}\\
				=&{}|\tau(a_j'e_i')-\frac{1}{n}\tau(a_j')|+\frac{2\varepsilon}{3}.
			\end{align*}
			For $\tau\in \mathrm{T}(A\rtimes_\alpha G)$, we use $\bar{\tau}$ to denote the restriction on $B$. So we can get a tracial state $\hat{\tau}\colon=\frac{\bar{\tau}}{\|\bar{\tau}\|}$ on $B$.
			We know that $\|\bar{\tau}\|=d_\tau(d)=\lim_{n\to\infty}\tau(d^\frac{1}{n})$. Since $d$ is positive contraction, we have
			\begin{equation}\label{eq10}
			\|\bar{\tau}\|=d_\tau(d)\leq 1.
			\end{equation}
			Therefore, by (\ref{eq10}) and (\ref{eq11}), for every $\tau\in\mathrm{T}(A\rtimes_\alpha G)$, we have
			\begin{align*}
				|\tau(a_je_i)-\frac{1}{n}\tau(a_j)|\leq&{}|\tau(a_j'e_i')-\frac{1}{n}\tau(a_j')|+\frac{2\varepsilon}{3}\\
				=&{}\|\bar{\tau}\||\hat{\tau}(a_j'e_i')-\frac{1}{n}\hat{\tau}(a_j')|+\frac{2\varepsilon}{3}\\
				=&{}d_\tau(d)|\hat{\tau}(a_j'e_i')-\frac{1}{n}\hat{\tau}(a_j')|+\frac{2\varepsilon}{3}\\
				<&{}|\hat{\tau}(a_j'e_i')-\frac{1}{n}\hat{\tau}(a_j')|+\frac{2\varepsilon}{3}.
			\end{align*}
		    Thus, by (\ref{eq11}), we have
			\begin{align*}
				\sup_{\tau\in\mathrm{T}(A\rtimes_\alpha G)}|\tau(a_je_i)-\frac{1}{n}\tau(a_j)|
				\leq&{}\sup_{\hat{\tau}\in\mathrm{T}(B)}|\hat{\tau}(a_j'e_i')-\frac{1}{n}\hat{\tau}(a_j')|+\frac{2\varepsilon}{3}\\
				<&{}\frac{\delta'}{2}+\frac{2\varepsilon}{3}<\varepsilon.
			\end{align*}
			By Proposition \ref{UPG}, $A\rtimes_\alpha G$ has uniform proerty $\Gamma$.
		\end{proof}
	\end{thm}
	
	\begin{cor}
		Let $A$ be a unital separable simple infinite dimensional C*-algebra which has uniform property $\Gamma$. Let $\alpha\colon G\to \mathrm{Aut}(A)$ be an action of a finite group which has the weak tracial Rokhlin property. Then the fixed point algebra $A^\alpha$ has uniform property $\Gamma$.
		\begin{proof}
			As in the proof of Theorem \ref{FG}, we know that $A^\alpha\cong pM_m(A\rtimes_\alpha G)p$ for some $m\in\mathbb{N}$ and $p\in M_m(A\rtimes_\alpha G)$. By Lemma \ref{TM} and Lemma \ref{HS}, we know that $A^\alpha$ has uniform property $\Gamma$.
		\end{proof}
	\end{cor}
	
	The second theorem uses the different method. We don't directly consider the crossed product but the fixed point algebra. Since there is a sequence of asymptotic homomorphisms from the original algebra to the fixed point algebra, the structural property of the sequence algebra of the fixed point algebra can be studied by considering the natural homomorphism from the original algebra to the sequence algebra of the fixed point algebra, and then we can use some technique to get the structural property of the fixed point algebra. Since the crossed product is stably isomorphic to the fixed point algebra, we will get the structural property of the crossed product.
	\begin{thm}\label{CG}
		Let $A$ be a unital separable simple infinite dimensional C*-algebra which has uniform property $\Gamma$. Let $\alpha\colon G\to \mathrm{Aut}(A)$ be an action of a second-countable compact group which has the tracial Rokhlin property with comparison. Then the fixed point algebra $A^\alpha$ has uniform property $\Gamma$.
		\begin{proof}
			Since $A$ has uniform property $\Gamma$, we know that  $\mathrm{T}(A)$ is nonempty. For $\tau\in\mathrm{T}(A)$, we can restrict it on $A^\alpha$ as a trace. Since $1_A=1_{A^\alpha}$, we know that $\tau\in\mathrm{T}(A^\alpha)$. Thus $\mathrm{T}(A^\alpha)$ is nonempty. This together with the unitality of $A^\alpha$ implies that $\mathrm{T}(A^\alpha)$ is compact. By Proposition \ref{UPG}, we need to show that for any finite subset $F\subset A^\alpha$, any $\varepsilon>0$, and any $n\in\mathbb{N}$, there exist pairwise orthogonal positive contractions $e_1,\dots,e_n\in A^\alpha$ such that
			\[\|e_ia-ae_i\|<\varepsilon \quad \mathrm{and} \quad \sup_{\tau\in\mathrm{T}(A^\alpha)}|\tau(ae_i)-\frac{1}{n}\tau(a_i)|<\varepsilon,
			\]
			for $i=1,\dots,n$ and $a\in F$.
			
			Without loss of generality, we may assume that $\|a\|\leq1$ for all $a\in F$. We choose $0<\delta<\frac{\varepsilon}{6}$. Since $A^\alpha$ is simple and not of type \uppercase\expandafter{\romannumeral1} (see \cite[Theorem 3.2]{MP2022} and \cite[Proposition 3.3]{MP2022} ), by \cite[Corollary 2.5]{NCP2016}, there exists a nonzero positive element $x\in A^\alpha$ such that $d_\tau(x)<\delta$ for all $\tau\in\mathrm{T}(A^\alpha)$.
			
			Since $A$ has uniform property $\Gamma$, by Proposition \ref{UPG}, with $n$, $F$, and $\delta$ as given, we have pairwise orthogonal positive contractions $e'_1,\dots,e'_n\in A$ such that
			\begin{equation}\label{eq5}
				\|e'_ia-ae'_i\|<\delta \quad \mathrm{and} \quad \sup_{\tau\in\mathrm{T}(A)}|\tau(ae'_i)-\frac{1}{n}\tau(a_i)|<\delta,
			\end{equation}
			for $i=1,\dots,n$ and $a\in F$.
			
			Set $F_1=\{e'_1,\dots,e'_n\}$. Suppose that $\{S_m\}$ is a sequence of finite sets in $A$ such that $\cup_{m\in\mathbb{N}}S_m$ is dense in $A$ and $F_1\subseteq S_m$ for all $m\in\mathbb{N}$. Apply Theorem \ref{MTF} for $\frac{1}{m}$, $x$, $S_m$, $F$, there exist a projection $p_m\in A^\alpha$ and a unital completely positive map $\varphi_m\colon A\to p_mA^\alpha p_m$ such that
			\begin{enumerate}
				\item $\varphi_m$ is an ($2$,$S_m$,$\frac{1}{m}$)-approximately multiplicative map.
				
				\item $\|p_ma-ap_m\|<\frac{1}{m}$ for all $a\in S_m$.
				
				\item $\|\varphi_m(a)-p_map_m\|<\frac{1}{m}$ for all $a\in F$.
				
				\item $1-p_m\precsim_{A^\alpha} (x-\frac{1}{2})_+$.
			\end{enumerate}
			Let $\pi_A\colon l^\infty(\mathbb{N},A^\alpha)\to (A^\alpha)_\infty$ be the quotient map. Define a homomorphism $\varphi\colon A\to q(A^\alpha)_\infty q$ by $\varphi(a)=\pi_A(\{\varphi_1(a),\varphi_2(a),\dots\})$ for all $a\in A$. Denote $p=(p_m)_{m=1}^\infty\in l^\infty(\mathbb{N},A^\alpha)$ and $q=\pi_A(p)$. We have
			\begin{enumerate}
				\setcounter{enumi}{4}
				\item $\varphi(1)=q$.
				
				\item $aq=qa$ for all $a\in A$.\label{item2.6}
				
				\item $\varphi(a)=qaq$ for all $a\in F$.\label{item2.7}
				
				\item $1-q\precsim_{(A^\alpha)_\infty}x$.\label{item2.8} 
			\end{enumerate}
			Then, for every $\tau\in\mathrm{T}(A^\alpha)$, it induces a tracial state on $(A^\alpha)_\infty$, we still use $\tau$ to denote it. Using (\ref{item2.6}) at the first step, (\ref{item2.7}) at the second step, (\ref{item2.8}) at the third step, we have
			\begin{align}
				&{}|\tau(a\varphi(e_i'))-\frac{1}{n}\tau(a)| \nonumber\\
				\leq&{}|\tau(qaq\varphi(e_i'))-\frac{1}{n}\tau(qaq)|+|\tau((1-q)a(1-q)\varphi(e_i'))-\frac{1}{n}\tau((1-q)a(1-q))| \nonumber\\
				\leq&{}|\tau(\varphi(ae_i'))-\frac{1}{n}\tau(\varphi(a))|+|0-\frac{1}{n}\tau(1-q)| \nonumber\\
				\leq&{}|\tau(\varphi(ae_i'))-\frac{1}{n}\tau(\varphi(a))|+|\frac{1}{n}d_\tau(x)| \nonumber\\
				<&{}|\tau(\varphi(ae_i'))-\frac{1}{n}\tau(\varphi(a))|+\delta\label{eq3}.
			\end{align}
			For $\frac{\varepsilon}{3}>0$ and $n\in\mathbb{N}$, we choose $\delta'<\frac{\varepsilon}{24}$ sufficiently small such that satisfying Lemma \ref{RO}. Then we can choose $m_0$ sufficiently large to get a unital completely positive map $\varphi_{m_0}\colon A\to p_{m_0}A^\alpha p_{m_0}$ such that
			\begin{enumerate}
				\setcounter{enumi}{8}
				\item $\varphi_{m_0}$ is an ($2$,$F_1$,$\delta'$)-approximately multiplicative map.\label{item2.9}
				
				\item $\|p_{m_0}a-ap_{m_0}\|<\delta'$ for all $a\in F_1$.
				
				\item $\|\varphi_{m_0}(a)-p_{m_0}ap_{m_0}\|<\delta'$ for all $a\in F$.\label{item2.11}
				
				\item $1-p_{m_0}\precsim_{A^\alpha} x$.
			\end{enumerate}
			and for all $\tau\in \mathrm{T}(A^\alpha)$,
			\begin{equation}\label{eq2}
				|\tau(a\varphi_{m_0}(e_i'))-\frac{1}{n}\tau(a)|\leq|\tau(a\varphi(e_i'))-\frac{1}{n}\tau(a)|+\delta'.
			\end{equation}
			By (\ref{item2.9}), we have $\|\varphi_{m_0}(e_i')\varphi_{m_0}(e_j')\|<\delta'$, for $i\neq j$. By Lemma \ref{RO}, there exist pairwise orthogonal positive contractions $e_1,\dots,e_n\in p_{m_0}A^\alpha p_{m_0}$ such that 
			\begin{equation}\label{eq6}
				\|e_i-\varphi_{m_0}(e_i')\|<\frac{\varepsilon}{3},
			\end{equation}
			for $i=1,\dots,n$.
			
			Thus for all $\tau\in \mathrm{T}(A^\alpha)$, we have
			\begin{equation}\label{eq1}
				|\tau(ae_i)-\frac{1}{n}\tau(a)|<|\tau(a\varphi_{m_0}(e_i'))-\frac{1}{n}\tau(a)|+\frac{\varepsilon}{3}.
			\end{equation}
			Since $\varphi$ is a unital homomorphism from $A\to q(A^\alpha)_\infty q$, for all $\tau\in\mathrm{T}(A^\alpha)$, we have $\frac{1}{\tau(q)}\tau\circ\varphi\in\mathrm{T}(A)$ and $\tau(q)\leq 1$. Thus, we have
			\begin{equation}\label{eq4}
				\sup_{\tau\in\mathrm{T}(A^\alpha)}|\tau(\varphi(ae_i'))-\frac{1}{n}\tau(\varphi(a))|\leq\sup_{\tau\in\mathrm{T}(A)}\tau(q)|\tau(ae_i')-\frac{1}{n}\tau(a)|\leq\sup_{\tau\in\mathrm{T}(A)}|\tau(ae_i')-\frac{1}{n}\tau(a)|.
			\end{equation}
			Therefore, using (\ref{eq1}) at the first step, (\ref{eq2}) at the second step, (\ref{eq3}) at the third step, (\ref{eq4}) at the fourth step, (\ref{eq5}) at the fifth step, we have
			\begin{align*}
				\sup_{\tau\in\mathrm{T}(A^\alpha)}|\tau(ae_i)-\frac{1}{n}\tau(a)|\leq&{}\sup_{\tau\in\mathrm{T}(A^\alpha)}|\tau(a\varphi_{m_0}(e_i'))-\frac{1}{n}\tau(a)|+\frac{\varepsilon}{3}\\
				\leq&{}\sup_{\tau\in\mathrm{T}(A^\alpha)}|\tau(a\varphi(e_i'))-\frac{1}{n}\tau(a)|+\delta'+\frac{\varepsilon}{3}\\
				<&{}\sup_{\tau\in\mathrm{T}(A^\alpha)}|\tau(\varphi(ae_i'))-\frac{1}{n}\tau(\varphi(a))|+\delta+\delta'+\frac{\varepsilon}{3}\\
				\leq&{}\sup_{\tau\in\mathrm{T}(A)}|\tau(ae_i')-\frac{1}{n}\tau(a)|+\delta+\delta'+\frac{\varepsilon}{3}\\
				<&{}2\delta+\delta'+\frac{\varepsilon}{3}\\
				<&{}\frac{\varepsilon}{3}+\frac{\varepsilon}{24}+\frac{\varepsilon}{3}<\varepsilon.
			\end{align*}
			For all $a\in F$, using (\ref{item2.11}) and $e_i\in  p_{m_0}A^\alpha p_{m_0}$, we have
			\begin{align*}
				\|ae_i-\varphi_{m_0}(a)e_i\|\leq\|p_{m_0}ap_{m_0}e_i-\varphi_{m_0}(a)e_i\|+\|(1-p_{m_0})a(1-p_{m_0})e_i\|<\delta',
			\end{align*}
			and, using (\ref{eq6}) and (\ref{item2.9}), we have
			\begin{align*}
				\|\varphi_{m_0}(a)e_i-\varphi_{m_0}(ae'_i)\|\leq&{}\|\varphi_{m_0}(a)e_i-\varphi_{m_0}(a)\varphi_{m_0}(e'_i)\|+\|\varphi_{m_0}(a)\varphi_{m_0}(e'_i)-\varphi_{m_0}(ae'_i)\|\\
				<&{}\frac{\varepsilon}{3}+\delta'.
			\end{align*}
			Thus we have
			\[\|ae_i-\varphi_{m_0}(ae'_i)\|\leq\|ae_i-\varphi_{m_0}(a)e_i\|+\|\varphi_{m_0}(a)e_i-\varphi_{m_0}(ae'_i)\|<\frac{\varepsilon}{3}+2\delta'.\]
			With the same argument, we have
			\[\|\varphi_{m_0}(e'_ia)-e_ia\|<\frac{\varepsilon}{3}+2\delta'.\]
			Since $\|e'_ia-ae'_i\|<\delta$, we have
			\begin{align*}
				\|ae_i-e_ia\|\leq&{}\|ae_i-\varphi_{m_0}(ae'_i)\|+\|\varphi_{m_0}(ae'_i)-\varphi_{m_0}(e'_ia)\|+\|\varphi_{m_0}(e'_ia)-e_ia\|\\
				<&{}\frac{2\varepsilon}{3}+4\delta'+\delta\\
				<&{}\frac{2\varepsilon}{3}+\frac{\varepsilon}{6}+\frac{\varepsilon}{6}=\varepsilon.
			\end{align*}
			By Proposition \ref{UPG}, $A^\alpha$ has uniform proerty $\Gamma$.
		\end{proof}
	\end{thm}
	
	\begin{cor}
		Let $A$ be a unital separable simple infinite dimensional C*-algebra which has stablized property $\Gamma$. Let $\alpha\colon G\to \mathrm{Aut}(A)$ be an action of a second-countable compact group which has the tracial Rokhlin property with comparison. Then the crossed product $A\rtimes_\alpha G$ has uniform property $\Gamma$.
		\begin{proof}
			Since $A^\alpha$ has uniform proerty $\Gamma$ and \cite[Porposition 2.6]{CE2021}, $A^\alpha$ has stabilised proerty $\Gamma$. Since $A^\alpha$ and $A\rtimes_\alpha G$ are stably isomorphic, by \cite[Theorem 2.10]{CE2021}, $A\rtimes_\alpha G$ has stablized property $\Gamma$. Again by \cite[Porposition 2.6]{CE2021}, $A\rtimes_\alpha G$ has uniform property $\Gamma$.
		\end{proof}
	\end{cor}

	\subsection*{Acknowledgments}
	The authors would like to thank the referees for their helpful comments.

\end{document}